\def\version{August 7, 2005}

\documentclass[reqno,11pt]{amsart}

\usepackage{amsmath}
\usepackage{amssymb}

\newfam\Bbbfam
\font\tenBbb=msbm10
\font\sevenBbb=msbm7
\font\fiveBbb=msbm5
\textfont\Bbbfam=\tenBbb
\scriptfont\Bbbfam=\sevenBbb
\scriptscriptfont\Bbbfam=\fiveBbb

\newcommand{\R}     {\mathbb{R}}
\newcommand{\Z}     {\mathbb{Z}}
\newcommand{\N}     {\mathbb{N}}
\renewcommand{\P}   {\mathbb{P}}

\newcommand{\E}     {\mathbb{E}}

\def\1{{\mathchoice {1\mskip-4mu\mathrm l} 
{1\mskip-4mu\mathrm l}
{1\mskip-4.5mu\mathrm l} {1\mskip-5mu\mathrm l}}}

\def\comment#1{}
\newtheoremstyle{thm}{2ex}{2ex}{\itshape\rmfamily}{}
{\bfseries\rmfamily}{}{1.7ex}{}

\newtheoremstyle{rem}{1.3ex}{1.3ex}{\rmfamily}{}
{\itshape\rmfamily}{}{1.5ex}{}

\newenvironment{Proof}[1]
{\vskip0.1cm\noindent{\bf #1. }}{\vspace{0.15cm}}

\newtheorem{theorem}{Theorem}[section]
\newtheorem{lemma}[theorem]{Lemma}
\newtheorem{prop}[theorem] {Proposition}

\theoremstyle{definition}

\newtheorem{step}{STEP}

\newcommand{\en}       {\end{equation}}
\newcommand{\eq}       {\begin{equation}}

\newcommand{\eqry}   {\begin{eqnarray}}
\newcommand{\enqry}   {\end{eqnarray}}
\newcommand{\eqarray}   {\begin{eqnarray}}
\newcommand{\enarray}   {\end{eqnarray}}
\newcommand{\eqarraystar} {\begin{eqnarray*}}
\newcommand{\enarraystar} {\end{eqnarray*}}
\newcommand{\bel}{\begin{lemma}}
\newcommand{\el}{\end{lemma}}
\newcommand{\bes}{\begin{step}}
\newcommand{\es}{\end{step}}
\newcommand{\bea}{\begin{array}}
\newcommand{\ea}{\end{array}}
\newcommand{\bpr}{\begin{proof}}
\newcommand{\epr}{\end{proof}}

\renewcommand{\section}{\secdef\sct\sect}
\newcommand{\sct}[2][default]{\refstepcounter{section}
\vspace{0.5cm}
\setcounter{equation}{0}
\centerline{ 
\scshape \arabic{section}.\ #1}
\vspace{0.3cm}}
\newcommand{\sect}[1]{
\vspace{0.5cm}
\centerline{\scshape\Large\bf #1}
\vspace{0.3cm}}

\renewcommand{\subsection}{\secdef \subsct\sbsect}
\newcommand{\subsct}[2][default]{\refstepcounter{subsection}
\nopagebreak
\vspace{0.5\baselineskip}
{\flushleft\bf \arabic{section}.\arabic{subsection}~\bf #1  }
\nopagebreak}
\newcommand{\sbsect}[1]{\vspace{0.1cm}\noindent
{\bf #1}\vspace{0.1cm}}

\renewcommand{\subsubsection}{%
\secdef \subsubsect\sbsbsect}
\newcommand{\subsubsect}[2][default]{%
\refstepcounter{subsubsection}
\nopagebreak
\vspace{0.1\baselineskip}
\nopagebreak
{\flushleft
\sffamily\slshape
\arabic{section}.\arabic{subsection}.\arabic{subsubsection}
\ %
\sffamily #1\/.}\ }
\newcommand{\sbsbsect}[1]{\vspace{0.1cm}\noindent
{\bf #1}\ }

\renewcommand{\d}{{\rm d}}

\newcommand{\eps}{\varepsilon}

\newcommand{\ssup}[1] {{\scriptscriptstyle{({#1}})}}

\newcommand{\smfrac}[2]{\textstyle{\frac {#1}{#2}}}

\newcommand{\Ocal}   {{\mathcal O }}

\begin{document}

\title[Deviations of a RWRS with stretched exponential tails]{\large Deviations of a random walk in a random
scenery with stretched exponential tails}

\author[Remco van der Hofstad, Nina Gantert and Wolfgang
        K{\"o}nig]{}
\maketitle

\thispagestyle{empty}
\vspace{0.2cm}

\centerline {{\sc Remco van der Hofstad$^1$\/, Nina Gantert$^2$ and
Wolfgang  K{\"o}nig$^3$ }}
\vspace{0.8cm}

\centerline {\em $^1$Department of Mathematics and Computer Science,}
\centerline{\em Eindhoven University of Technology,}
\centerline
{\em Post Box 513, 5600 MB Eindhoven, The
Netherlands.}
\centerline{\em
    {\tt r.w.v.d.hofstad@tue.nl }}
\vspace{0.3cm}
\centerline {\em $^2$Institut f\"ur Mathematische Statistik, Universit\"at M\"unster}
\centerline {\em Einsteinstrasse 62, D-48149 M\"unster, Germany}
\centerline{\em
    {\tt gantert@math.uni-muenster.de}}
\vspace{0.3cm}
\centerline{\em $^3$Mathematisches Institut, Universit\"at Leipzig,}
\centerline{\em Augustusplatz 10/11, D-04109 Leipzig, Germany}
\centerline{\tt koenig@math.uni-leipzig.de}

\vspace{0.5cm}

\centerline{\small(\version)}
\vspace{0.5cm}

\begin{quote}
  {\small {\bf Abstract:}} Let $(Z_n)_{n\in\N_0}$ be a
  $d$-dimensional {\it random walk  in  random scenery}, i.e.,
  $Z_n=\sum_{k=0}^{n-1}Y_{S_k}$ with $(S_k)_{k\in\N_0}$ a
  random walk in $\Z^d$ and $(Y_{z})_{z\in\Z^d}$ an i.i.d.~scenery,
  independent of the walk.
  We assume that the random variables $Y_{z}$ have a stretched exponential tail. In particular, they do not possess exponential moments.
  We identify the speed and the rate of the logarithmic decay of
  $\P(\frac 1n Z_n>t_n)$ for all sequences $(t_n)_{n\in\N}$ satisfying a certain lower bound. This complements results of \cite{GKS04}, where it was assumed that $Y_{z}$ has exponential moments of all orders. In contrast to the situation \cite{GKS04}, the event $\{\frac 1n Z_n>t_n\}$ is not realized by a homogeneous behavior of the walk's local times and the scenery, but by many visits of the  walker to a particular site and a large value of the scenery at that site. This reflects a well-known extreme behavior typical for random variables having no exponential moments.
\end{quote}



\setcounter{section}{0}
\section{Introduction}
 \label{Intro}

\subsection{The model}

\noindent Let $S=(S_n)_{n\in\N_0}$ be a random walk on $\Z^d$
starting at the origin (more precisely, $S=(S_n)_{n\in\N_0}$ is a sequence of partial sums of i.i.d.~$\Z^d$-valued random variables). Defined on the same probability space,
let $Y=(Y_{z})_{z\in\Z^d}$ be an i.i.d.~sequence of  random variables,
independent of the walk. We refer to $Y$ as the {\it random scenery.}
Then the process $(Z_n)_{n\in\N}$ defined by
$$
Z_n=\sum_{k=0}^{n-1} Y_{S_k},\qquad n\in\N,
$$
where $\N=\{1, 2, \ldots\},$ is called a {\it random walk in random scenery} (RWRS), sometimes
also referred to as the {\it Kesten-Spitzer random walk in random scenery}, see \cite{KS79}. An
interpretation is as follows. If a random walker pays
$Y_{z}$ units at any time he/she visits the site $z$, then $Z_n$
is the total amount he/she pays by time $n-1$.
We denote by $\P$ the underlying
probability measure
and by $\E$ the corresponding expectation.

The random walk in random scenery has been introduced and
analyzed for dimension $d\not=2$ by H.~Kesten and F.~Spitzer \cite{KS79} and
by E.~Bolthausen \cite{B89} for $d= 2$. Under the assumptions that the walk is in the domain of attraction of Brownian motion and that
$Y_{0}$ has expectation zero and variance $\sigma^2\in(0,\infty)$, their results imply that
\begin{equation}
    \label{weakconv}
   \frac 1n Z_n\approx a_n=
    \begin{cases}n^{-\frac 14}&\mbox{if }d=1,\\
    (\frac n{\log n})^{-\frac 12}&\mbox{if }d=2,\\
    n^{-\frac 12}&\mbox{if }d\ge 3.
    \end{cases}
\end{equation}
More precisely,
$\frac 1{na_n}Z_n$ converges in distribution towards some non-degenerate
random variable. The limit is Gaussian in $d\ge 2$ and a
convex combination of Gaussians (but not Gaussian) in $d=1$. This
can be roughly explained as follows.
In terms of the so-called {\it local times\/} of the walk and its {\it range},
\begin{equation}
    \label{loctim}
    \ell_n(z)=\sum_{k=0}^{n-1} \1_{\{S_k=z\}},\quad R_n=\{S_0,S_1,\dots,S_{n-1}\},\qquad n\in\N, \;\;
    z\in\Z^d,
\end{equation}

\noindent the random walk in random scenery may be identified as
\begin{equation}
    \label{Znrepr}
    Z_n=\sum_{z\in R_n}Y_{z}\ell_n(z).
\end{equation}

\noindent Hence, conditionally on the
random walk, $Z_n$ is, for dimension $d\ge 3$, a sum of
$\Ocal(n)$ independent copies of  finite multiples of $Y_{0}$,
and hence it is plausible that $n^{-1/2}Z_n$  converges to a
normal variable.  The same assertion with logarithmic
corrections is also plausible in $d=2$. However,  in $d=1$, $Z_n$ is
roughly a sum of
$\Ocal(n^{1/2})$ copies of independent variables with
variances of order $\Ocal( n)$, and this suggests the
normalization in \eqref{weakconv} as well as a non-Gaussian limit.

In this paper, we  analyse deviations $\{\frac 1n Z_n>t_n\}$ for sequences
$(t_n)_n$ of positive numbers satisfying $t_n\gg a_n$, by which we mean that $\lim_{n\to\infty}t_n/a_n=\infty$. The problem of deviations of the random walk in random scenery and also of the continuous version, Brownian motion in a random scenery on $\R^d$, have gained interest in recent years. One reason is that the interplay between the trajectory and the medium displays a rich behavior and is therefore mathematically appealing. Furthermore, the continuous version of this problem appears in the asymptotic  analysis of diffusions in a Gaussian shear flow drift (see \cite{CP01}, \cite{Ca04}, \cite{AC03} and the references cited therein). Furthermore, there is a tight methodological relationship to the parabolic Anderson model (see the survey in \cite{GK05}), where one studies the asymptotics of the exponential moments of the random walk in random scenery, in the continuous setting or in the spatially discrete, but time-continuous setting.
In fact, precise logarithmic asymptotics for the decay of the probability of $\{\frac 1n Z_n>t_n\}$ correspond to moment asymptotics of the parabolic Anderson model with suitably rescaled scenery.

The main question is the description
of the \lq optimal\rq\ behavior of the walk and of the scenery to meet the
event $\{\frac 1n Z_n>t_n\}$ in the \lq cheapest\rq\ way.
So far, only random sceneries
having exponential moments of all orders have been considered
In \cite{AC02},  the random sceneries are bounded, and in \cite{GKS04},
the random sceneries
have exponential moments of all orders.
In these cases, if the tail of the scenery decays fast enough w.r.t.~ the dimension, it turns out that the optimal behavior is homogeneous in the sense that, in a certain centered ball with $n$-dependent radius,
all the walker's local times and all the scenery values grow unboundedly, each with its
appropriate speed. The exponential decay rate of the probability
of $\{\frac 1n Z_n>t_n\}$ is characterized in terms of a variational problem.
If the tail of the random sceneries decays slower (but still having exponential moments of all orders), the optimal strategy is different, we refer to \cite{AC05} for recent results.

\subsection{Our main result}

\noindent In the present paper, we study the deviation problem in the
case where the scenery has a {\it stretched exponential tail}. In particular, it does not have any positive exponential moments. It is known
that the cheapest way for a sum of i.i.d.~stretched-exponential random variables to attain a huge value is to make just one of these variables as huge as required, and the others do not contribute. Our main result shows that a similar picture appears for the random walk in random scenery.

We turn to a description of the results of this paper. Our assumptions on the random i.i.d.~scenery $(Y_z)_{z\in\Z^d}$ are the following:

\medskip

\noindent{\bf Centering Assumption.} The random variable $Y_0$ satisfies
    \begin{equation}
    \E [Y_{0}] =0,  \quad \E [Y_{0}^2] = \sigma^2 < \infty,
    \end{equation}
and

\medskip

\noindent{\bf Tail Assumption.} {\it There is a constant $q\in (0,1)$ and a
slowly varying function $D\colon (0,\infty)\to(0,\infty)$ such that
    \begin{equation}\label{tailass}
    \log \P (Y_{0}> t)\sim -D(t) t^q,\qquad \mbox{as }t\to\infty.
    \end{equation}
Moreover, the map $t\mapsto D(t) t^{q-1}$ is eventually decreasing, and}
    \begin{equation}\label{gammadef}
    \gamma(a)=\lim_{t\rightarrow \infty} \frac{D(t^{a+o(1)})}{D(t)} \in (0,\infty) \qquad\mbox{\it {exists for every} }a\in (0,1).
    \end{equation}
    \medskip

(We write $b_t\sim c_t$ for $t\to\infty$ if $\lim_{t\to\infty}b_t/c_t=1$.)
In fact, \eqref{gammadef} implies that $D$ is slowly varying.
Consequently, $\gamma(a)$ is a power of $a$, but we are not going to use this fact. The Tail Assumption says that the upper tails of the scenery variables are of Weibull type, modulo some technical regularity assumption, and they have no positive exponential moment.

Our assumptions on the random walk are the following. For $d\leq 2$, the walk is
assumed to be  recurrent. In $d=2$, we furthermore assume that $\sup_{k\in\N}k\P(S_k=0)<\infty$.
Furthermore, we assume that the limits
\begin{equation}\label{Cddef}
K_d=\begin{cases}
\lim\limits_{n\to\infty}n^{-\frac 12}\E[\ell_n(0)]&\mbox{if }d=1,\\
\lim\limits_{n\to\infty}\frac 1{\log n}\E[\ell_n(0)]&\mbox{if }d=2,
    \end{cases}
    \end{equation}
exist in $(0,\infty)$. This includes the case of simple random walk with $K_1=2\pi^{-1/2}$ and $K_2=1/\pi$.

Let
\begin{equation}\label{betandef}
 \beta_n(t)=\begin{cases}n^{\frac q{q+2}}t^{\frac {2q}{q+2}}\left(D(nt^2)\gamma(\frac 1{q+2})\right)^{\frac{1}{q+2}}
&\mbox{if }d=1,\\
    \big((q+1)nt\big)^{\smfrac{q}{q+1}}(\log\frac n{t^q})^{-\smfrac{q}{q+1}}\left(D(nt)\gamma(\frac 1{q+1})\right)^{\frac 1{q+1}}&\mbox{if }d=2,\\
(nt)^{\smfrac{q}{q+1}}\left(D(nt)\gamma(\frac 1{q+1})\right)^{\frac{1}{q+1}}&\mbox{if }d\geq 3.
    \end{cases}
    \end{equation}

Our main result is the following:
\begin{theorem}
\label{main}
Fix $d\geq 1$ and a sequence $(t_n)_{n\in\N}$ of positive numbers such that
\begin{equation}\label{tnass}
t_n\geq n^{-r}\qquad\mbox{for some }r<\begin{cases} \frac {1-q}{4-q}&\mbox{if }d=1,\\
\frac {1-q}{2}&\mbox{if }d\geq 2.
\end{cases}
\end{equation}
Then, as $n\to\infty$,
\begin{equation}\label{onedim}
 \log \P(Z_n> nt_n)\sim -\beta_n(t_n)\times\begin{cases}
    \big(4K_1^2/q\big)^{\frac {2q}{q+2}}(2+q)&\mbox{if }d=1,\\
    (K_2/q)^{\frac {q}{q+1}}(1+q)&\mbox{if }d= 2.
    \\
    (-\frac 1q\log f_0)^{\frac {q}{q+1}}(1+q)&\mbox{if }d\geq 3,
    \end{cases}
\end{equation}
where $f_0=\P(S_n = 0\mbox{ for some }n\in\N)$ is the return probability of the random walk.
\end{theorem}

Note that in \eqref{tnass} only a lower bound on $t_n$ is imposed. Our assumptions on $t_n$  leave a gap to the scale $a_n$ of the limit law in \eqref{weakconv}.
We think that the result persists to a wider range of $t_n$'s, but not to sequences $t_n$ that are too close to $a_n$. For more detailed comments, we refer to Section \ref{Open}.

\subsection{Outline of the proof}\label{sec-outline}

\noindent An explanation of Theorem~\ref{main} and of its proof is as follows. Recall that stretched exponential random variables have the characteristic property that a sum of $n$ independent copies has the same large deviation behavior as just one of them. That is, for i.i.d.~random variables $Y_1,Y_2,Y_3,\dots$ having the same distribution as our scenery variables, we have, for $t>0$ fixed,
\begin{equation}
\label{iidas}
\log \P \left(\sum\limits_{i=1}^n Y_i > nt\right) \sim \log\P(Y_0>nt)\sim -(nt)^{q}D(n),\qquad n\to\infty.
\end{equation}
This is proved in \cite{Na69}; the (critical) upper bound in \eqref{iidas} is also a consequence of Lemma \ref{weightas} below.
For the random walk in random scenery with stretched exponential tails, it turns out in our first result that the large deviation behavior of $Z_n=\sum_z Y_z\ell_n(z)$ is also governed by just one summand:

\begin{prop}\label{onesite}
Under the Centering Assumption and the Tail Assumption, for any sequence $(t_n)_{n\in\N}$ satisfying \eqref{tnass}, and for any $\eps>0$,
\begin{equation}\label{only0up}
\limsup_{n\to\infty}\frac 1{\beta_n(t_n)}\log\frac{\P(Z_n>nt_n)}{\P(Y_0\ell_n(0)>nt_n(1-\eps))}\leq 0
\end{equation}
and
\begin{equation}\label{only0lo}
\liminf_{n\to\infty}\frac 1{\beta_n(t_n)}\log\frac{\P(Z_n>nt_n)}{\P(Y_0\ell_n(0)>nt_n(1+\eps))}\geq 0\, .
\end{equation}
\end{prop}

Hence, it suffices to identify the large deviation behaviors of $\frac 1{nt_n} Y_0$ and of  $\ell_n(0)$ and to combine the two in an appropriate manner. For doing this, it is convenient to introduce a new scale function $1\ll \alpha_n\ll nt_n$ and to look at large deviation principles for $\frac {\alpha_n}{nt_n}Y_0$ and $\frac 1{\alpha_n}\ell_n(0)$. It is clear from the Tail Assumption that $\frac {\alpha_n}{nt_n}Y_0$ satisfies a large deviation principle on $(0,\infty)$ with rate function $y\mapsto y^q$ and speed $(nt_n/\alpha_n)^qD(nt_n/\alpha_n)$, i.e., as $n\to\infty$,
\begin{equation}\label{YLDP}
\log\P(\smfrac {\alpha_n}{nt_n}Y_0>y)\sim -y^q\Big(\frac{nt_n}{\alpha_n}\Big)^qD\big(\smfrac{nt_n}{\alpha_n}\big),\qquad y>0.
\end{equation}
Furthermore, the moderate deviations for the local time $\ell_n(0)$ are identified as follows.

\begin{lemma}\label{LoctimLDP}
Let $\alpha_n\ll n$ for any dimension  $d\in\N$, and $\sqrt n\ll\alpha_n$ in $d=1$, $\log n\ll\alpha_n$ in $d=2$, and $1\ll\alpha_n$ in $d\geq 3$. Then, as $n\to\infty$,
\begin{equation}\label{LDPelln(0)}
\log\P(\ell_n(0)>\alpha_n)\sim -\begin{cases}K_1^2\smfrac {\alpha_n^2}n&\mbox{if }d=1,\\
K_2\smfrac{\alpha_n}{\log\smfrac n{\alpha_n}}&\mbox{if }d=2,\\
-(\log f_0) \alpha_n&\mbox{if }d\geq 3,\end{cases}
\end{equation}
where $K_d$ is defined in \eqref{Cddef}, and $f_0=\P(S_n=0\mbox{ for some }n\in\N)$ is the return probability.
\end{lemma}

It remains to pick $\alpha_n$ such that the two speeds in \eqref{YLDP} and \eqref{LDPelln(0)} coincide, i.e., such that
\begin{equation}\label{speedequal}
\left(\frac {nt_n}{\alpha_n}\right)^qD\left(\smfrac {nt_n}{\alpha_n}\right)\sim\begin{cases}\smfrac {\alpha_n^2}n&\mbox{if }d=1,\\
\smfrac{\alpha_n}{\log\smfrac n{\alpha_n}}&\mbox{if }d=2,\\
\alpha_n&\mbox{if }d\geq 3.\end{cases}
\end{equation}
This is guaranteed by the choice
    \begin{equation}
    \label{alphadef}
    \alpha_n
    =\begin{cases}
    n^{\smfrac{q+1}{q+2}}t_n^{\smfrac{q}{q+2}}\left(D(nt_n^2)\gamma(\frac 1{q+2})\right)^{\smfrac{1}{q+2}}
    &\mbox{if }d=1,\\
    (nt_n)^{\smfrac{q}{q+1}}\left(\frac{\gamma(\frac 1{q+1}) D(nt_n)\log \frac n{t_n^q}}
   {q+1}\right)^{\smfrac{1}{q+1}}&\mbox{if }d=2,\\
    (nt_n)^{\smfrac{q}{q+1}} \left(D(nt_n)\gamma(\frac 1{q+1})\right)^{\smfrac{1}{q+1}}&\mbox{if }d\geq 3,
    \end{cases}
    \end{equation}
where we have also used \eqref{gammadef}.  The speeds of the two principles  in \eqref{YLDP} and \eqref{LDPelln(0)} are then both equal to the speed $\beta_n(t_n)$ in \eqref{betandef}. It remains to combine the two principles for $Y_0$ and $\ell_n(0)$, which is elementary. This ends the explanation of the proof of Theorem~\ref{main}. We see that the event $\{Z_n>nt_n\}$ is optimally met by sceneries having $Y_0$ of order $nt_n/\alpha_n$ and random walks having $\ell_n(0)$ of order $\alpha_n$ with $\alpha_n$ in \eqref{alphadef}.


    The proof of Proposition~\ref{onesite} is in  Section~\ref{sec-proofs}, and the proof Lemma~\ref{LoctimLDP} and the completion of the proof of Theorem~\ref{main} are in Section~\ref{sec-finish}. In Section \ref{Open}, we give some open questions and conjectures.

\section{Approximation of $Z_n$ by $Y_0\ell_n(0)$}\label{sec-proofs}

\noindent In Section~\ref{sec-PropProof} we prove Proposition~\ref{onesite}. As an important pre-step, we give a generalization of (\ref{iidas})
for weighted sums of random variables in Section~\ref{sec-conditioned}.

\subsection{A conditional estimate}\label{sec-conditioned}

\noindent The following lemma can be seen as a conditional upper estimate for random walk in random scenery, given the random walk.

\begin{lemma}
\label{weightas}
Assume $(Y_i)_{i\in\N}$ is a sequence of i.i.d.~random variables satisfying the Centering Assumption and the Tail Assumption.
Fix a sequence $(t_n)_n$ of positive numbers satisfying \eqref{tnass}, and abbreviate $m_n=n t_n^{(2-q)/(1-q)}$.
Then, for any  $\eta>0$, any sufficiently large $n$, every  $r\in\{1,\dots, n\}$ and any choice of $l_1,\dots,l_r\in[1,\infty)$ satisfying $\sum_{i=1}^r l_i = n$ and $L\equiv\max_{i=1,\dots,r}l_i\leq [nt_n\wedge m_n]^{1-\eta}$,
\begin{equation}
\label{tailweight}
\P \left(\sum\limits_{i=1}^r l_i Y_i > nt_n\right) \leq
\exp\left(-\left(\frac {nt_n}L\right)^qD\left(\smfrac {nt_n}L\right)(1 -4\varepsilon)\right).
\end{equation}
\end{lemma}

Remark: the proof given below only uses the fact that $D$ is slowly varying and not the (stronger)
 property \eqref{gammadef}.

\begin{Proof}{Proof}
We begin with
\begin{equation}
\label{splitmax}
\P \left(\sum\limits_{i=1}^r l_i Y_i > nt_n\right) \leq
\P \left(L\max\limits_{1 \leq i \leq r} Y_i > nt_n\right)
+\P \left(\sum\limits_{i=1}^r l_i Y_i > nt_n, \max\limits_{1 \leq i \leq r}(l_i Y_i) \leq nt_n\right)\, .
\end{equation}
With the help of \eqref{tailass}, the first term on the r.h.s.\ of (\ref{splitmax}) can, for all large $n$ and all $r\in\{1,\dots,n\}$, be estimated by
\begin{equation}
\label{splitmax1}
\P \left(L\max\limits_{1 \leq i \leq r} Y_i > nt_n\right)
\leq  n\exp\left(-\left(\frac {nt_n}L\right)^q D(\smfrac {nt_n}L)(1-\eps)\right)\, .
\end{equation}
For estimating the second term on the r.h.s.\ of (\ref{splitmax}), we use the Markov inequality. For any $\lambda >0$ (to be determined later),
\begin{equation}
\label{splitmax2}
\begin{aligned}
\P \left(\sum\limits_{i=1}^r l_i Y_i > nt_n, \max\limits_{1 \leq i \leq r}(l_i Y_i) \leq nt_n\right)
&\leq {\rm e}^{-\lambda nt_n}\prod\limits_{i=1}^r \E \left[{\rm e}^{\lambda l_i Y_i}\1_{\{l_iY_i \leq nt_n\}}\right]\\
&={\rm e}^{-\lambda nt_n}\prod\limits_{i=1}^r \big[\Lambda_i^{\ssup{1}}(n) +
\Lambda_i^{\ssup{2}}(n)\big],
\end{aligned}
\end{equation}
where
\begin{equation}
\Lambda_i^{\ssup{1}}(n) = \E \left[{\rm e }^{\lambda l_i Y_i}\1_{\{l_iY_i < \lambda^{-1}\}}\right]\qquad\mbox{and}\qquad
\Lambda_i^{\ssup{2}}(n) = \E \left[{\rm e}^{\lambda l_i Y_i} \1_{\{\lambda^{-1} \leq l_iY_i \leq nt_n\}}\right].
\end{equation}
Fix $i\in\{1,\dots,r\}$. We have to estimate $\Lambda_i^{\ssup{1}}(n)$ and  $\Lambda_i^{\ssup{2}}(n)$.  Using first the inequality ${\rm e}^u \leq 1+u+u^2$ for $u<1$ and then $1+u \leq {\rm e}^u$, and taking into account that $\E [Y_i ] =0$ and $\E [Y_i^2 ] = \sigma^2$, we have
\begin{equation}
\label{Lambda1est}
\Lambda_i^{\ssup{1}}(n) \leq 1 + \lambda^2l_i^2\E [Y_i^2 ] \leq {\rm e}^{\lambda^2l_i^2\sigma^2}\, .
\end{equation}
To estimate $\Lambda_i^{\ssup{2}}(n)$, we use the following estimate, which is valid for any random variable $X$ and any $\lambda >0$ and $0 < T_1 < T_2 < \infty$,
\begin{equation}
\label{intparts}
\E \left[{\rm e}^{\lambda X} \1_{\{T_1 \leq X \leq T_2\}}\right] \leq \int\limits_{T_1}^{T_2}\lambda
{\rm e}^{\lambda s}\P ( X> s)\,\d s + {\rm e}^{\lambda T_1}\P (X \geq  T_1) .
\end{equation}
Hence,
\begin{equation}
\Lambda_i^{\ssup{2}}(n) \leq \lambda\int\limits_{\lambda^{-1}}^{nt_n}{\rm e}^{\lambda s} \P(l_iY_i >s)\,\d s
+ {\rm e} \,\P (l_iY_i \geq\lambda^{-1}).
\end{equation}
We now determine $\lambda=\lambda_n$ by
    \begin{equation}
    \label{lambdandef}
    \lambda_n =  \frac 1{nt_n}\left(\frac{nt_n}L\right)^q D(\smfrac{nt_n}L)(1-2\varepsilon).
    \end{equation}
Note that $\lim_{n\to\infty}\lambda_n=0$. Recalling our assumption in \eqref{tailass}, we obtain, for all large $n$, all $r\in\{1,\dots,n\}$ and all $i\in\{1,\dots,r\}$,
\begin{equation}
\label{Lambda2est}
\Lambda_i^{\ssup{2}}(n) \leq \lambda_n \int\limits_{\lambda_n^{-1}}^{nt_n} {\rm e}^{\lambda_n s}{\rm e}^{-D(s/l_i)(1-\eps) s^q l_i^{-q}} \,\d s
+{\rm e}^{1-D(\lambda_n^{-1}l_i^{-1})(1-\eps)\lambda_n^{-q}l_i^{-q}}.
\end{equation}
We are going to estimate the integral on the right hand side of \eqref{Lambda2est}. We claim that, for any large $n\in\N$ and all $1\leq i\leq r\leq n$,
    \begin{equation}
    \label{claim}
    \lambda_n s - D(\smfrac{s}{l_i}) (1-\eps)\left(\smfrac s{ l_i}\right)^{q} \leq -\varepsilon D(\smfrac{s}{l_i}) \left(\smfrac s{ l_i}\right)^{q} \qquad \hbox{for any } s
    \in [\lambda_n^{-1}, nt_n].
    \end{equation}
Define $f(s)=D(s) s^{q-1}$, then the claim in \eqref{claim} is equivalent to
    \eq
    (1-2\varepsilon)f(s/l_i)\geq\lambda_n l_i,\qquad\hbox{for any } s
    \in [\lambda_n^{-1}, nt_n].
    \en
We note that
    \eq
    \inf_{s\in [\lambda_n^{-1}, nt_n]}s/l_i\geq (\lambda_n l_i)^{-1} \geq (\lambda_n L)^{-1}
    =(1-2\eps)^{-1}f(\smfrac{nt_n}{L})^{-1}\rightarrow \infty\qquad\mbox{as }n\to\infty,
    \en
since $\frac{nt_n}{L}\rightarrow \infty$ by the assumption that $L\leq (nt_n)^{1-\eta}$,
and $f(s)\rightarrow 0$ as $s\rightarrow \infty$. Recall that $f$ is eventually decreasing by our Tail Assumption. Hence, $s\mapsto f(s/l_i)$ is decreasing in $[\lambda_n^{-1}, nt_n]$ for all sufficiently large $n$.
Therefore, to prove the claim, it is enough to verify $\lambda_n l_i-(1-2\varepsilon)f(s/l_i)
\leq 0$ only for the right end-point, $s=nt_n$.
For this, we note that
$$
\lambda_n l_i-(1-2\varepsilon)f({\smfrac{nt_n}{l_i}}) \leq  (1-2\varepsilon)[f({\smfrac{nt_n}{L}})
-f({\smfrac{nt_n}{l_i}})] \leq 0,
$$
again by monotonicity of $f$. This proves the claim in \eqref{claim}.

We pick some $\widetilde q\in(0,q)$. Hence, for $n$ large enough, we obtain, using \eqref{claim}, the substitution $u=s/l_i$, and the estimate $D(u)u^q \geq u^{\widetilde q}$ for large $u$, (the latter follows since $D$ is slowly varying)
\begin{equation}
\begin{aligned}
\lambda_n \int\limits_{\lambda_n^{-1}}^{nt_n}{\rm e}^{\lambda_n s}{\rm e}^{-D(s/l_i)(1-\eps) s^q l_i^{-q}} \,\d s
&\leq
\lambda_n l_i\int\limits_{(\lambda_nl_i)^{-1}}^{nt_n/l_i} {\rm e}^{-\varepsilon D(u) u^q } \,\d u\leq\lambda_n l_i\int\limits_{(\lambda_nl_i)^{-1}}^{\infty}{\rm e}^{-\varepsilon u^{\widetilde{q} }} \,\d u\\
&= \frac{\lambda_n l_i}{\widetilde{q}}\int\limits_{(\lambda_nl_i)^{-\widetilde{q}}}^{\infty} t^{\widetilde{q}^{-1} -1}{\rm e}^{-\varepsilon t} \,\d t\leq {\rm e}^{-\smfrac12 \varepsilon (\lambda_nl_i)^{-\widetilde{q}}}.
\end{aligned}
\end{equation}
Going back to (\ref{Lambda2est}), we have, using that $\lambda_n \rightarrow 0$ and
$\lambda_n l_i\leq \lambda_n L\to 0$,
    \begin{equation}
    \Lambda_i^{\ssup{2}}(n) \leq  {\rm e}^{-\smfrac12 \varepsilon (\lambda_nl_i)^{-\widetilde{q}}}
    + {\rm e}^{1-D(\lambda_n^{-1}l_i^{-1})(1-\eps)(\lambda_nl_i)^{-q}} = o((\lambda_nl_i)^2)
    \qquad \hbox{ for } n \to \infty,
    \end{equation}
uniformly in $i$.
Hence, for $n$ large enough,
\begin{equation}
\Lambda_i^{\ssup{1}}(n)+\Lambda_i^{\ssup{2}}(n)
\leq {\rm e}^{\lambda_n^2l_i^2\sigma^2}+\varepsilon \lambda_n^2l_i^2\sigma^2
\leq  {\rm e}^{(1+\varepsilon)\lambda_n^2l_i^2\sigma^2}.
\end{equation}
Therefore, using the inequality $\sum_{i=1}^r l_i^2 \leq \sum_{i=1}^r l_i L=nL$ and recalling the choice of $\lambda_n$ in \eqref{lambdandef},
\begin{equation}
\label{lambdaprodest}
\begin{aligned}
\prod\limits_{i=1}^r \big[\Lambda_i^{\ssup{1}}(n)+\Lambda_i^{\ssup{2}}(n)\big]
&\leq
\exp\left((1+\varepsilon)\lambda_n^2\sigma^2\sum\limits_{i=1}^r l_i^2\right) \\
&\leq \exp\left((1+\varepsilon)\sigma^2 n^{2q-1}t_n^{2q-2}L^{-2q+1}
D(\smfrac{nt_n}L)^2(1-2\varepsilon)^2\right)\\
&\leq \exp\left(2\sigma^2\left(\frac {nt_n}L\right)^qD(\smfrac{nt_n}L)\frac 1{t_n}\left(\frac {nt_n}L\right)^{q-1}D(\smfrac{nt_n}L)\right)\\
&=\exp\Big(2\sigma^2\left(\frac {nt_n}L\right)^qD(\smfrac{nt_n}L)\Big(\frac{L}{m_n}\Big)^{1-q}D(\smfrac{nt_n}L)\Big),
\end{aligned}
\end{equation}
where we recall that $m_n=n t_n^{(2-q)/(1-q)}$. By our assumption that $L\leq m_n^{1-\eta}$, and since $D$ is slowly varying, the right hand side can be estimated,
for all large $n$, against $\exp(\eps \left(\frac {nt_n}L\right)^qD(\smfrac{nt_n}L))$.
Using this in \eqref{splitmax2} and recalling \eqref{lambdandef}, we obtain that
\begin{equation}
\P \left(\sum\limits_{i=1}^r l_i Y_i > nt_n, \max\limits_{1 \leq i \leq r}(l_i Y_i) \leq nt_n\right)
\leq \exp\left(- \left(\frac {nt_n}L\right)^q(1-3\eps)D(\smfrac{nt_n}L)\right)
\end{equation}
Together with \eqref{splitmax1} and \eqref{splitmax}, we arrive at the assertion.
\qed
\end{Proof}

%

\subsection{Proof of Proposition~\ref{onesite}}\label{sec-PropProof}

\noindent We begin with \eqref{only0up}. Pick $\alpha_n$ as in \eqref{alphadef}. We again use the abbreviation $m_n=n t_n^{(2-q)/(1-q)}$. Denote by $L_n=\max_{z\in\Z^d}\ell_n(z)$ the maximal local time of the random walk. Fix a small $\eta>0$. Estimate
\begin{equation}\label{onesite1}
\P(Z_n>nt_n)\leq \P(L_n>m_n^{1-\eta})+\P(Z_n>nt_n, L_n\leq m_n^{1-\eta}).
\end{equation}
Observe that
$$
\P(L_n>m_n^{1-\eta})\leq \P(\ell_n(z)>m_n^{1-\eta}\mbox{ for some }z\in \Z^d,|z|\leq n)\leq (2n+1)^d \P(\ell_n(0)>m_n^{1-\eta}),
$$
since $x\mapsto \ell_n(x)$ is stochastically maximal in $x=0$. We now choose $\eta>0$ so small that $m_n^{1-\eta}\gg\alpha_n$. This is possible because of \eqref{tnass}. Then, with the help of Lemma~\ref{LoctimLDP} and \eqref{alphadef}, we see that the first term in \eqref{onesite1} is negligible:
$$
\limsup_{n\to\infty}\frac 1{\beta_n(t_n)}\log\P(L_n>m_n^{1-\eta})=-\infty.
$$
In order to treat the second term in \eqref{onesite1}, we apply Lemma~\ref{weightas} to $Z_n$,
recalling \eqref{Znrepr}, and condition on the local times of the random walk (recall \eqref{loctim}). Fix $\eps>0$ so small that $(1-5\eps/q)^q<1-4\eps$. We condition on $\ell_n(\cdot)$ and obtain from Lemma~\ref{weightas}, for all large $n$, on the event $\{ L_n\leq m_n^{1-\eta}\}$,
\begin{equation}
\P(Z_n>nt_n\mid \ell_n)= \P \left(\sum\limits_{x \in R_n}Y_x \ell_n(x) > nt_n \,\Big|\, \ell_n\right)\leq \exp\left(-\left(\smfrac{nt_n}{L_n}\right)^q D(\smfrac{nt_n}{L_n})(1-4\varepsilon)\right).
\end{equation}
Using again the Tail Assumption, we obtain, for all large $n$,
\begin{equation}
\P(Z_n>nt_n\,\mid \,\ell_n)\leq  \P \left(Y_0 L_n > nt_n (1-5\varepsilon/q)\,\big|\, \ell_n\right).
\end{equation}
Integrating over $\ell_n$ on the event $\{ L_n\leq m_n^{1-\eta}\}$, we conclude that
\begin{equation}
\label{annupper}
\begin{aligned}
\P(Z_n>nt_n, L_n\leq m_n^{1-\eta})&\leq \sum_{x\in R_n}\P\left(Y_0 L_n> nt_n (1-5\varepsilon/q), L_n=\ell_n(x),L_n\leq m_n^{1-\eta} \right)\\
&\leq (2n+1)^d\P\left(Y_0 \ell_n (0) > nt_n (1-5\varepsilon/q)\right).
\end{aligned}
\end{equation}
This implies \eqref{only0up}.

We turn now to the proof of \eqref{only0lo}. Abbreviate $\widetilde{Z_n}= \sum_{x \in R_n\setminus\{0\}}Y_x \ell_n(x)$ and pick some $\eps\in(0,1/(2\sigma^2))$. Then we have
\begin{equation}\label{annupper1}
\begin{aligned}
\P(Z_n>nt_n)&\geq \P\left(Y_0 \ell_n(0)> nt_n(1+\varepsilon), \widetilde{Z_n}>-\varepsilon nt_n\right)\\
&\geq \E\left[\P\left(Y_0 \ell_n(0)> nt_n(1+\varepsilon)\,\big|\,\ell_n\right)\1_{\{L_n\leq \varepsilon^3nt_n^2\}}\P\big(\widetilde{Z_n}>-\varepsilon nt_n\,\big|\,\ell_n\big)\right].
\end{aligned}
\end{equation}
Using the Chebyshev inequality, we estimate the last term as follows.
    \begin{equation}
    \begin{aligned}
    \P\left(\widetilde{Z_n} > -\varepsilon nt_n\,\big|\,\ell_n\right)
    &=1-\P\left(\widetilde{Z_n} \leq -\varepsilon nt_n\,\big|\,\ell_n\right)
    \geq
    1-\frac{1}{(\varepsilon nt_n)^2}\hbox{ Var}(\widetilde{Z_n}\,\big|\,\ell_n)\\
    &=1-\frac{1}{(\varepsilon nt_n)^2}\sum_{z\in R_n\setminus\{0\}}\ell_n(z)^2\sigma^2
    \geq 1-\frac{\sigma^2}{\eps^2}\,\frac {L_n}{nt_n^2}.
    \end{aligned}
    \end{equation}
Hence, on $\{L_n\leq \varepsilon^3 nt_n^2\}$, we have $\P(\widetilde{Z_n} > -\varepsilon nt_n\,\big|\,\ell_n)\geq \frac 12$ for all sufficiently large $n$. This gives in \eqref{annupper1}
\begin{equation}\label{annupper2}
\begin{aligned}
\P(Z_n>nt_n)&\geq \frac 12 \P\left(L_n\leq \varepsilon^3 nt_n^2,Y_0 \ell_n(0)> nt_n(1+\varepsilon)\right)\\
&\geq \frac 12 \Big[\P(Y_0 \ell_n(0)> nt_n(1+\varepsilon))-\P\left(L_n> \varepsilon^3 nt_n^2\right)\Big].
\end{aligned}
\end{equation}
We estimate $\P(L_n> \varepsilon^3 nt_n^2)\leq n^d\P(\ell_n(0)> \varepsilon^3 nt_n^2)$. In our proof of Theorem~\ref{main} in Section~\ref{sec-ThmProof} below we will see that $\P(Y_0 \ell_n(0)> nt_n(1+\varepsilon))\geq {\rm e}^{-\Ocal(\beta_n(t_n))}$. Observe that $(nt_n^2)^{1-\eta}\gg \alpha_n$ for some $\eta>0$. Indeed, this holds as soon as $t_n\geq n^{-r}$ with $r<(4+q)^{-1}$ in $d=1$ and $r<(2+q)^{-1}$
in $d\geq 2$, and this is implied by \eqref{tnass}. Therefore, Lemma~\ref{LoctimLDP} implies that $\P(L_n> \varepsilon^3 nt_n^2)$ is much smaller than $\P(Y_0 \ell_n(0)> nt_n(1+\varepsilon))$.
Hence, the last line of \eqref{annupper2} can be estimated from below by $\frac 14 \P(Y_0 \ell_n(0)> nt_n(1+\varepsilon))$, and this completes the proof of \eqref{only0lo}.
\qed


\section{Moderate deviations for the local time, and the proof of Theorem~\ref{main}}\label{sec-finish}

\noindent We prove the moderate deviations statement for the local time $\ell_n(0)$ (Lemma~\ref{LoctimLDP}) in Section~\ref{sec-Loctim}, and we complete the proof of Theorem~\ref{main} in Section~\ref{sec-ThmProof}.

\subsection{Proof of Lemma~\ref{LoctimLDP}}\label{sec-Loctim}

\noindent The statement \eqref{LDPelln(0)} in $d=1$ follows from \cite[Theorem~2]{C01} with $f(X_k) = \1_{\{X_k =0\}}$, $a(n) \sim K_1 \sqrt{n}$, $p = \smfrac12$.

In $d=2$, it follows from \cite[Theorem~1]{GZ98}. In the notation in \cite{GZ98}, $g(n) = \E[\ell_n(0)]\sim K_2\log n$ and $\alpha_n=\psi(n) g(n)$. We note that in \cite{GZ98}, it is assumed that $n\mapsto \psi(n)$
is non-decreasing. However, an inspection of the proof shows that the monotonicity is not used at all, but only that $\psi(n)$ does not vanish as $n\to\infty$.

In $d\geq 3$, the proof of \eqref{LDPelln(0)} is easily done as follows. Let  $T_0 = 0<T_1<T_2< \dots$ denote the subsequent times at which the walker hits the origin, i.e., $T_{i} = \inf\{n> T_{i-1}\colon S_n =0\}$ for $i\in\N$. Then $f_0=\P(T_1<\infty)$, and we have
\begin{equation}
\label{returns}
\P\left(\ell_n(0)>  \alpha_n\right) \leq \P\left(\lim_{m\to\infty}\ell_m(0)>\alpha_n\right)=\P(T_{\alpha_n}<\infty)=f_0^{\alpha_n},
\end{equation}
which is the upper bound in \eqref{LDPelln(0)}. To prove the lower bound, note that, for $n\to\infty$,
\begin{equation}
\begin{aligned}
\P\left(\ell_n(0) > \alpha_n\right)&\geq \P\left(T_i-T_{i-1}<\smfrac{n}{\alpha_n}\,\forall i=1,\dots,\alpha_n\right)=\P(T_1< \smfrac{n}{\alpha_n})^{\alpha_n}
=\left(f_0-o(1)\right)^{\alpha_n}\\
&=f_0^{\alpha_n}{\rm e}^{o(\alpha_n)},
\end{aligned}
\end{equation}
since $\smfrac{n}{\alpha_n}\rightarrow \infty$. This completes the
proof of Lemma~\ref{LoctimLDP}.
\qed

\subsection{Proof of Theorem~\ref{main}}\label{sec-ThmProof}

\noindent From now on, we pick $\alpha_n$ as in \eqref{alphadef} and $\beta_n(t_n)$ as in \eqref{betandef} with $t=t_n$. Recall that \eqref{speedequal} is satisfied, and note that $\beta_n(t_n)$ is given by \eqref{betandef}. From \eqref{YLDP}, we in particular have, as $n\to\infty$,
\begin{equation}\label{YLDPbeta}
\log\P\left(\frac {\alpha_n}{nt_n} Y_0>y\right)\sim -y^q\beta_n(t_n),\qquad y>0.
\end{equation}
Replacing $\alpha_n$ in \eqref{LDPelln(0)} by $x\alpha_n$ for some $x>0$, we obtain the large deviation statement
\begin{equation}\label{locandY1}
\lim_{n\to\infty}\frac 1{\beta_n(t_n)}\log\P\left(\frac{1}{\alpha_n}\ell_n(0)> x\right) = -I_\ell(x),\qquad x>0,
\end{equation}
where $I_\ell(x)=K_1 x^2$ in $d=1$, $I_\ell(x)=K_2 x$ in $d= 2$ (recall \eqref{Cddef}) and $I_\ell(x)=-x\log f_0$ in $d\geq 3$ (recall that $f_0$ is the return probability).

%
%

The large deviation principles in \eqref{YLDPbeta} and \eqref{locandY1}, together with  \cite[Ex.~4.2.7]{DZ93}, imply that the distributions of $\smfrac1{nt_n} Y_{0}\ell_n(0)$ satisfy a large deviation principle on $(0,\infty)$ with speed $\beta_n(t_n)$ and rate function
$$
\widetilde I(s)= \inf\limits_{y,x\in(0,\infty)\colon yx > s} [y^q+I_\ell(x)],
$$
that is,
\begin{equation}\label{LDPYell}
\lim_{n\rightarrow \infty}\frac{1}{\beta_n(t_n)}\log\P\left(Y_{0}\ell_n(0)> snt_n\right) = -\widetilde I(s),\qquad s>0.
\end{equation}
Therefore, it remains to determine $\widetilde I(1)$.
It is not hard to see that $\widetilde I(1)$ is equal to the constant on the right hand side of \eqref{onedim}. Hence, Proposition~\ref{onesite} completes the proof of Theorem \ref{main}.
\qed

%
%
%
\section{Heuristics for small deviations}
\label{Open}
Let us discuss the necessity of our assumption in \eqref{tnass}, which leaves
a gap to the scale $a_n$ of the limit law in \eqref{weakconv}. We believe that
our main result in \eqref{onedim} persists to a wider range of $t_n$'s, but not
all the way down to $a_n$. The main reason is that
one way to realize the event $\{Z_n>n t_n\}$ is to let the random walk 
behave like free random walk,
while the scenery variables on the range of the walk are all of order $t_n$. As it turns out, for $t_n$ sufficiently close
to $a_n$, this strategy yields a lower bound on $\P(Z_n>n t_n)$
that is larger than the asymptotics in Theorem \ref{main}. In particular,
we see that Lemma~\ref{weightas} (which is an important ingredient
of the proof of the upper bound of \eqref{onedim}) breaks down in
this regime.

Let us explain this more closely, first in the case $d\geq 3$. If $n^{-\frac 12}\ll
t_n \ll n^{-(1-q)/(2-q)}$, in contrast to Lemma~\ref{weightas} with $r=n$ and
$L=1$, a sum of i.i.d.\ random variables with tails given by (\ref{tailass}) satisfies
a moderate deviation principle of central limit theorem (CLT) type \cite{EL03},
that is,
    \begin{equation}\label{moddev}
    \lim_{n \to\infty}\frac{1}{nt_n^2} \log \P \left(\sum\limits_{i=1}^n  Y_i > nt_n\right) =
    -\frac{1}{2}\,.
    \end{equation}
Hence, we obtain a lower bound for $\log \P(Z_n>nt_n)$ of order $nt_n^2$
by requiring that the walk's range is of order $n$ (this has probability $e^{-\Ocal(1)}$)
and that the scenery
performs a CLT type moderate deviation on the vertices in the range.
Further restricting $t_n$ to satisfy $n^{-\frac 12}\ll t_n \ll
n^{-1/(q+2)}$, we have found a cheaper strategy than the one
of Theorem~\ref{main}, since $nt_n^2\ll (nt_n)^{q/(q+1)}$.
This shows that the asymptotics in \eqref{onedim} does not hold
for all $a_n\ll t_n\leq n^{-r}$ with $r\leq\frac 1{q+2}$ (this upper bound on $r$
is smaller than the lower bound on $r$ in \eqref{tnass}).
We expect that for $d=2$,
the same argument applies apart from logarithmic corrections.

In one dimension, the situation is slightly different. We obtain a
lower bound for $\P(Z_n>nt_n)$ by additionally requiring that the walk's range
and most of the local times in this range are of order $n^{\frac12}$.
The probability for this is again $e^{-\Ocal(1)}$.
Conditionally on this behavior of the walk, $Z_n$ is in distribution
roughly equal to $n^{1/2}\sum_{i=1}^{n^{1/2}} Y_i$.
Using \eqref{moddev}, we see that, for $n^{-\frac 14}\ll t_n
\ll n^{-\frac 12(1-q)/(2-q)}$,  the conditional probability of
$\{Z_n > nt_n\}$ is not smaller than $\exp\{-\Ocal(n^{1/2}t_n^2)\}$.
Further restricting to $n^{-\frac 14}\ll t_n\ll n^{-\frac18 (2-q)}$,
we have found a cheaper strategy than the one of Theorem~\ref{main}.
Indeed, the exponential speed in Theorem~\ref{main} is $n^{q/(q+2)}t_n^{2q/(q+2)}$, which is much larger than the speed $n^{1/2}t_n^2$ we obtained above.
This shows that the asymptotics in \eqref{onedim} does not hold for all
$a_n\ll t_n\leq n^{-r}$ with $r\leq\frac 18(2-q)$ (this upper bound on $r$
is again smaller than the lower bound on $r$ in \eqref{tnass}).
\vskip0.5cm

\subsection*{Acknowledgements.} This work began during a visit of NG to EURANDOM.
The work of RvdH was supported in part by The Netherlands
Organisation for Scientific Research (NWO). WK gratefully
acknowledges the Heisenberg grant awarded by the German Science Foundation (DFG).

\end{document}